\documentclass[runningheads]{lncse}

\usepackage{amsmath}
\usepackage{amsfonts}
\usepackage{epsfig}
\usepackage{graphics}

\usepackage{bm}
\usepackage{subfig}
\usepackage{tikz}
\usetikzlibrary{plotmarks}
\usetikzlibrary{shapes,arrows}
\usetikzlibrary{positioning}
\usetikzlibrary{trees,mindmap}
\usetikzlibrary{fadings,shapes.arrows,shadows}
\usetikzlibrary{decorations.pathreplacing}
\usetikzlibrary{calc}
\usepackage{pgfplots}

\begin{document}

\title{Numerical investigation on the fixed-stress splitting scheme for Biot's equations:
Optimality of the tuning parameter}

\titlerunning{Optimality of the fixed-stress splitting parameter}

\author{ Jakub W. Both\inst{1} \and Uwe K\"ocher\inst{2}   }

\authorrunning{J.W. Both et al.}   

\institute{
University of Bergen, Bergen, Norway, {\tt jakub.both@uib.no}
\and Helmut-Schmidt-University, University of the Federal Armed Forces Hamburg, Germany, {\tt uwe.koecher@hsu-hh.de}
}

\maketitle

\begin{abstract}
We study the numerical solution of the quasi-static linear Biot's equations
solved iteratively by the fixed-stress splitting scheme. In each iteration
the mechanical and flow problems are decoupled, where the flow problem is
solved by keeping an artificial mean stress fixed. This introduces a numerical
tuning parameter which can be optimized. We investigate numerically the 
optimality of the parameter and compare our results with physically and 
mathematically motivated values from the literature, which commonly only
depend on mechanical material parameters. We demonstrate, that the optimal 
value of the tuning parameter is also affected by the boundary conditions
and material parameters associated to the fluid flow problem suggesting
the need for the integration of those in further mathematical analyses
optimizing the tuning parameter.
\end{abstract}

\section{Introduction}

The coupling of mechanical deformation and fluid flow in porous media is 
relevant in many applications ranging from environmental to biomedical engineering.
In this paper, we consider the simplest possible fully coupled model given by the 
quasi-static Biot's equations~\cite{both_mini9:biot1941general}, coupling classical and well-studied 
subproblems from linear elasticity and single phase flow in fully saturated porous media. 
Due to the complex structure of the coupled problem, the development of monolithic solvers
is not trivial and topic of current research. Hence, instead of developing
new simulation tools for the coupled problem, due to their simplicity and flexibility
splitting methods have been very attractive recently allowing the use of independent, tailored simulators for both
subproblems.
Among various iterative splitting schemes, one of the most prominent schemes is
the physically motivated fixed-stress splitting scheme~\cite{both_mini9:settari} based on solving sequentially the mechanics
and flow problems while keeping an artificial mean stress fixed in the latter.
\textcolor{black}{From an abstract point of view, the splitting scheme is a linearization scheme 
employing positive pressure stabilization, a concept also applied for the linearization of other problems
as, e.g., the Richards equation~\cite{both_mini9:list:2016}.}
Addressing the physical formulation, the definition of the artificial mean stress includes a user-defined tuning parameter.
It can be chosen a priori such that the resulting fixed-stress splitting scheme
is unconditionally stable in the sense of a von Neumann stability analysis~\cite{both_mini9:kim:2011b}
and it is globally contractive~\cite{both_mini9:bause:2017,both_mini9:bause:2017b,both_mini9:both:2017,both_mini9:mikelic:2013}.
Suggested values for the tuning parameter from literature are either physically motivated~\cite{both_mini9:kim:2011b} 
or mathematically motivated~\cite{both_mini9:bause:2017,both_mini9:bause:2017b,both_mini9:both:2017,both_mini9:mikelic:2013}.
The latter works prove theoretically global contraction of the scheme,
allowing to optimize the resulting theoretical contraction rate, and hence, proposing
a value for the tuning parameter with suggested, better performance than for the physically
motivated parameters.
In general, the suggested values for the tuning parameter given in the literature do not 
necessarily yield a minimal number of iterations, which for strongly coupled problems is crucial, 
as then the performance of the fixed-stress splitting scheme is very sensitive to the 
choice of the tuning parameter. 
We note that the mentioned tuning parameters depend solely on mechanical material
parameters. However, practically, it is known that the physical character of
the problem governed by boundary conditions also affects the performance of the 
scheme~\cite{both_mini9:kim:2011b}, introducing the main difficulty finding an 
optimal tuning parameter which would yield a minimal number of iterations.
We note that the fixed-stress splitting scheme can also be applied as a preconditioner
for Krylov subspace methods solving the Biot's equation in a monolithic fashion. 
In this case, performance is less sensitive with respect to the tuning parameter.

In this work, we investigate numerically whether the optimal tuning parameter
obtained by simple trial and error is closer related to the mathematically or
the physically motivated parameters.
Furthermore, we investigate whether the optimal tuning parameter is also dependent on 
more than only mechanical properties. For this purpose, we perform a numerical study 
enhancing a test case from~\cite{both_mini9:bause:2017} and measure performance of 
the fixed-stress splitting scheme for different tuning parameters. Our main results are:
\begin{itemize}
  \item[$\bullet$] Boundary conditions affect the optimality of the tuning parameter.
  
  \item[$\bullet$] Fluid flow parameters affect the optimality of the tuning parameter. 
  
  \item[$\bullet$] Both should be included in the mathematical analysis allowing
  to derive theoretically an optimal tuning parameter.
\end{itemize}

\section{Linear Biot's equations}
We consider the quasi-static Biot's equations~\cite{both_mini9:biot1941general},
modeling fluid flow in a deformable, linearly elastic porous medium $\Omega\subset\mathbb{R}^d$,
$d\in\{1,2,3\}$, fully saturated by a slightly compressible fluid. Using mechanical
displacement $\bm{u}$, fluid pressure $p$ and fluid flux $\bm{q}$ as primary variables,
on the space-time domain $\Omega\times(0,T)$,
the governing equations written in a three-field formulation read
\begin{align}
 -\bm{\nabla}\cdot\left(2\mu \bm{\varepsilon}(\bm{u}) + \lambda \bm{\nabla}\cdot \bm{u} \bm{I} - \alpha p \bm{I}\right) &= \rho_\mathrm{b} \bm{g}, \label{both_mini9:eq:biot:u} \\
 \partial_t \left( \frac{p}{M} + \alpha \bm{\nabla}\cdot \bm{u} \right) + \bm{\nabla}\cdot \bm{q} &= 0,                    \label{both_mini9:eq:biot:p} \\
 \frac{\eta}{k}\bm{q} + \bm{\nabla} p &= \rho_\mathrm{f} \bm{g}.                                         \label{both_mini9:eq:biot:q}
\end{align}
Eq.~\eqref{both_mini9:eq:biot:u} describes balance of momentum at each time, Eq.~\eqref{both_mini9:eq:biot:p}
describes mass conservation and Eq.~\eqref{both_mini9:eq:biot:q} describes Darcy's law. 
Here, $\bm{\varepsilon}(\bm{u})=\frac{1}{2}\left( \bm{\nabla} \bm{u} + \bm{\nabla} \bm{u}^\top\right)$ is the linearized strain tensor,
$\mu$, $\lambda$ are the Lam\'e parameters (equivalent to Young's modulus $E$ and Poisson's ratio $\nu$
via $\mu=\frac{E}{2(1+\nu)}$ and $\lambda=\frac{E\nu}{(1+\nu)(1-2\nu)}$), $\alpha$ is the Biot coefficient, $M$ is the Biot 
modulus, $\rho_\mathrm{f}$ is the fluid density, $\rho_\mathrm{b}$ is the bulk density, $k$ is the absolute permeability, 
$\eta$ the fluid viscosity and $\bm{g}$ is the gravity vector. In this work, we
assume isotropic, homogeneous materials, i.e., all material parameters are constants.

The system~\eqref{both_mini9:eq:biot:u}--\eqref{both_mini9:eq:biot:q} is closed by postulating
initial conditions $\bm{u}=\bm{u}_0$, $p=p_0$ on $\Omega\times\{0\}$, satisfying
Eq.~\eqref{both_mini9:eq:biot:u}, and boundary conditions $\bm{u}=\bm{u}_\mathrm{D}$ on $\Gamma_\mathrm{D,m}\times(0,T)$, 
$(2\mu \bm{\varepsilon}(\bm{u}) + \lambda \bm{\nabla}\cdot \bm{u} \bm{I} - \alpha p \bm{I})\cdot \bm{n}=\bm{\sigma}_\mathrm{N}$ on $\Gamma_\mathrm{N,m}\times(0,T)$,
$p=p_\mathrm{D}$ on $\Gamma_\mathrm{D,f}\times(0,T)$, $\bm{q}\cdot \bm{n}=q_\mathrm{N}$ on $\Gamma_\mathrm{N,f}\times(0,T)$
on partitions $\Gamma_\mathrm{D,\star}\cup\Gamma_\mathrm{N,\star}=\partial\Omega$, $\star\in\{m,f\}$,
where $\bm{n}$ is the outer normal on $\partial\Omega$.

Here and in the remaining paper, we omit introducing a corresponding variational formulation 
and suitable function spaces, as they appear naturally. For details, we refer to our 
works~\cite{both_mini9:bause:2017,both_mini9:both:2017}.

\section{Fixed-stress splitting scheme}
We solve the coupled Biot's equations~\eqref{both_mini9:eq:biot:u}--\eqref{both_mini9:eq:biot:q} iteratively
using the fixed-stress splitting scheme~\cite{both_mini9:settari}, which decouples
the mechanics and fluid flow problems.
Each iteration, defining the approximate solution $(\bm{u},p,\bm{q})^i$, $i\in\mathbb{N}$, consists of two steps. 
First, the flow problem is solved assuming a fixed artificial, volumetric stress $\sigma_v=K_\mathrm{dr}\bm{\nabla}\cdot\bm{u} - \alpha p$,
where $K_\mathrm{dr}$ is a tuning parameter, which will be discussed in the scope
of this paper: Given $(\bm{u},p,\bm{q})^{i-1}$, find $(p,\bm{q})^{i}$ satisfying
\begin{align}
 \left( \frac{1}{M} + \frac{\alpha^2}{K_\mathrm{dr}} \right)\partial_t p^i + \bm{\nabla}\cdot \bm{q}^i &= \frac{\alpha^2}{K_\mathrm{dr}}  \partial_t p^{i-1} - \alpha \partial_t \bm{\nabla}\cdot \bm{u}^{i-1},  \label{both_mini9:eq:fs:p}\\
 \frac{\eta}{k}\bm{q}^i + \bm{\nabla} p^i &= \rho_\mathrm{f} \bm{g},			                                            						                 \label{both_mini9:eq:fs:q}
\end{align}
including corresponding initial and boundary conditions. 
Second the mechanics problem is solved with updated flow fields: Find $\bm{u}^i$ satisfying
corresponding boundary conditions and
\begin{align}
 -\bm{\nabla}\cdot\left(2\mu \bm{\varepsilon}(\bm{u}^i) + \lambda \bm{\nabla}\cdot \bm{u}^i \bm{I} - \alpha p^i \bm{I}\right) &= \rho_\mathrm{b} \bm{g}. \label{both_mini9:eq:fs:u}
\end{align}

\subsection*{The tuning parameter $K_\mathrm{dr}$}\label{both_mini9:Sec:FS}

The fixed stress splitting scheme can be interpreted as a two block Gauss-Seidel
method with an educated predictor for mechanical displacement used in the solution
of the flow problem. More precisely, the mechanics problem~\eqref{both_mini9:eq:biot:u} is solved inexactly 
for the volumetric deformation by reduction to the one-dimensional equation 
\begin{align}\label{both_mini9:eq:fs:ansatz}
 K_\mathrm{dr}\bm{\nabla}\cdot(\bm{u}^i - \bm{u}^{i-1}) - \alpha (p^i-p^{i-1})=0
\end{align}
and inserted into the flow equation~\eqref{both_mini9:eq:biot:p}. 
In the special case of nearly incompressible materials, i.e., $\mu/\lambda\rightarrow 0$ or $\nu\rightarrow 0.5$,
Eq.~\eqref{both_mini9:eq:biot:u} yields $\bm{\nabla} \partial_t (\lambda \bm{\nabla}\cdot  \bm{u} - \alpha p) \approx \textbf{0}$.
Hence, we expect the ansatz~\eqref{both_mini9:eq:fs:ansatz} to be nearly exact for $K_\mathrm{dr}=\lambda$,
yielding a suitable tuning parameter for nearly incompressible materials.

An exact inversion of Eq.~\eqref{both_mini9:eq:biot:u} for the volumetric deformation would be given by the divergence of 
a Green's function and thus would be defined locally and depend on fluid pressure, 
geometry, material parameters and boundary conditions, both associated with the mechanical subproblem.
However, due to lack of a priori knowledge and simplicity, in the literature,
the considered inexact inversion includes only fluid pressure and mechanical material parameters
introduced via the tuning parameter $K_\mathrm{dr}$.
Selected values are $K_\mathrm{dr}^\text{1D}$, $K_\mathrm{dr}^\text{2D}$ $K_\mathrm{dr}^\text{3D}$,
cf.~\cite{both_mini9:kim:2011b}, $K_\mathrm{dr}^{2\times\lambda}$, cf.~\cite{both_mini9:bause:2017,both_mini9:bause:2017b,both_mini9:mikelic:2013}, and 
$K_\mathrm{dr}^{2\times d\text{D}}$, cf.~\cite{both_mini9:both:2017,both_mini9:mikelic:2013}, defined by
\begin{align*}
 K_\mathrm{dr}^{d^\star\text{D}}=\frac{2\mu}{d^\star}+\lambda,\ d^\star\in\{1,2,3\},\qquad
 K_\mathrm{dr}^{2\times\lambda}=2\lambda, \qquad
 K_\mathrm{dr}^{2\times d\text{D}}&=2 \, K_\mathrm{dr}^{d\text{D}}.
\end{align*}
The choice $K_\mathrm{dr}^{d^\star\text{D}}$ is purely physically motivated and equals 
the bulk modulus of a $d^\star$-dimensional material. Independent of the spatial dimension $d$,
for uniaxial compression, biaxial compression or general deformations, we choose $d^\star=1,2,3$, respectively.
Following~\cite{both_mini9:kim:2011b}, if not known better a priori, choose $d^\star=d$. 
Eq.~\eqref{both_mini9:eq:fs:ansatz} with 
$K_\mathrm{dr}=K_\mathrm{dr}^{d\text{D}}$ corresponds to fixing the trace of the physical,
poroelastic stress tensor. 
The choices $K_\mathrm{dr}^{2\times\lambda}$ and $K_\mathrm{dr}^{2\times d\text{D}}$ have resulted from 
optimization of the obtained theoretical contraction rate. Those analyses have in common
that global convergence is guaranteed for $0 \leq K_\mathrm{dr} \leq K_\mathrm{dr}^{2\times\lambda}$
and $0 \leq K_\mathrm{dr}\leq K_\mathrm{dr}^{2\times d\text{D}}$, hence, the latter also covers convergence
for the physical choices. Additionally, the analyses indicate that the larger 
$K_\mathrm{dr}$ the faster convergence, suggesting that the mathematically motivated parameters
should yield better performance than the physically motivated parameters. 
In the following section we investigate this statement numerically.

\section{Numerical study -- Optimal tuning parameter $K_\mathrm{dr}$}

\begin{figure}[ht]
\centering
\subfloat[Test cases 1a/b/c.]{\includegraphics[width=0.35\textwidth]{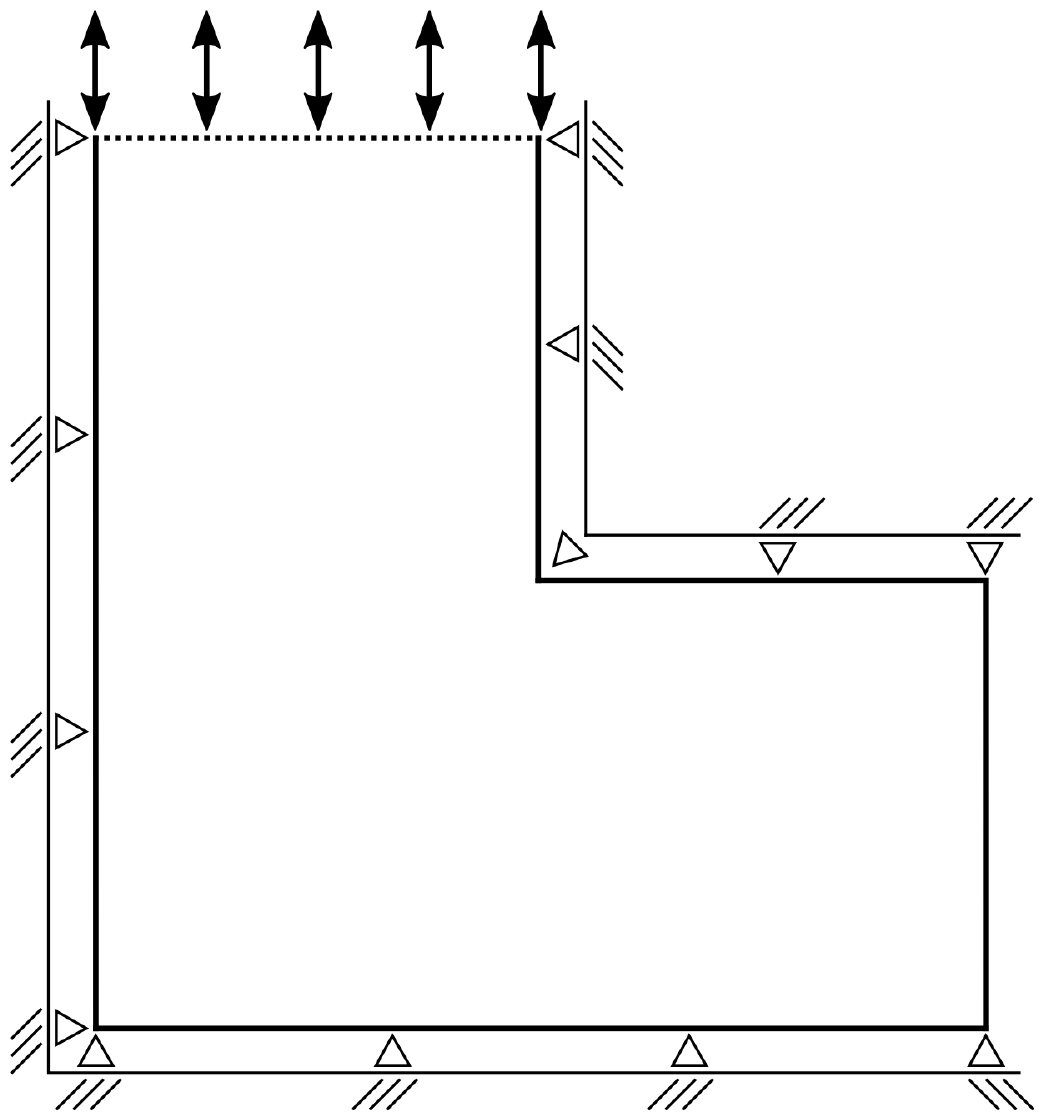}}
\hspace{1.5cm}
\subfloat[Test case 2.]{\includegraphics[width=0.35\textwidth]{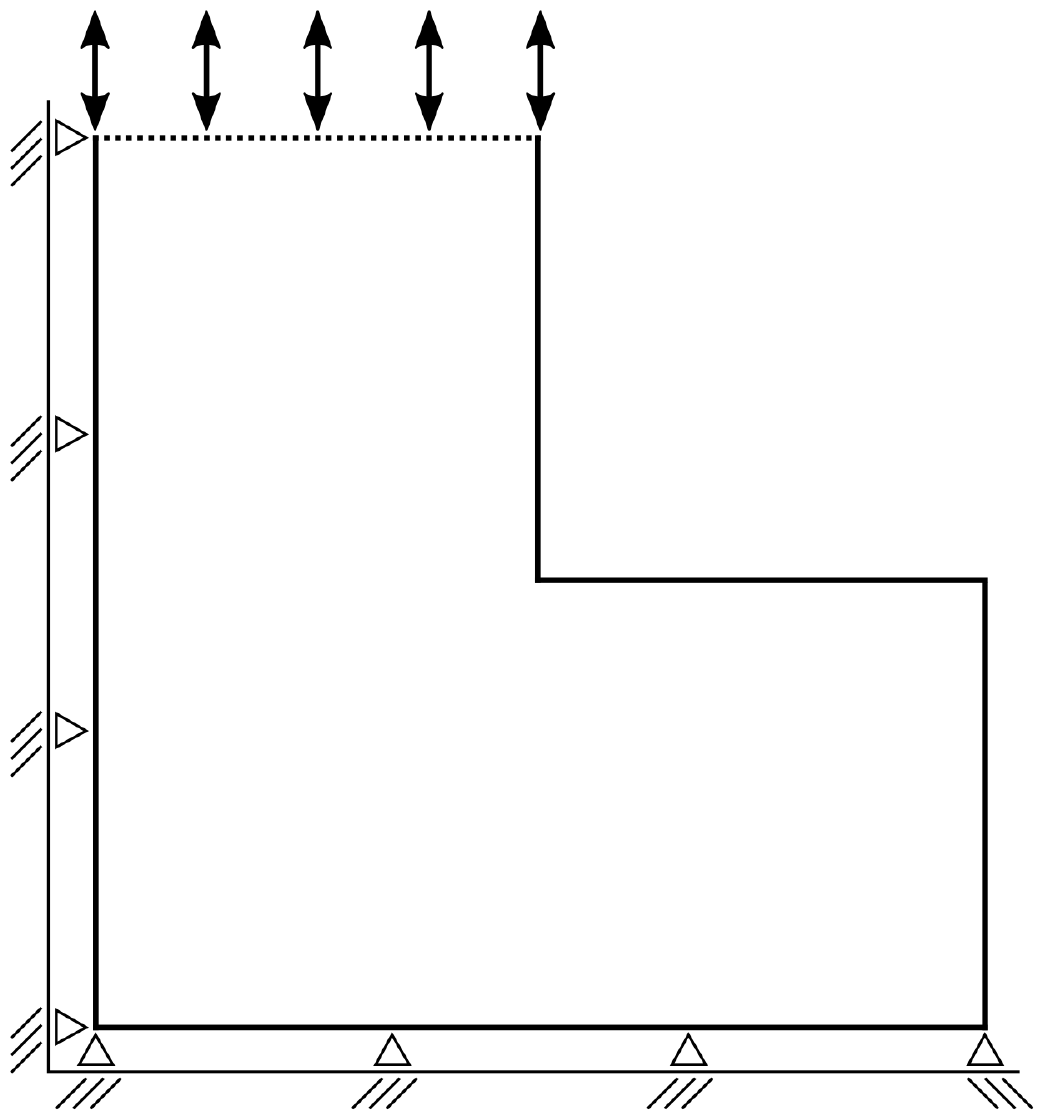}}
\caption{\label{both_mini9:figure:geometry} Geometry and boundary conditions employed in the numerical study.}
\end{figure}

We perform a numerical parameter study analyzing the optimality of the tuning parameter $K_\mathrm{dr}$
in the view of the performance of the fixed-stress splitting scheme measured in terms of number of 
fixed-stress iterations.
Inspired by a test setting from~\cite{both_mini9:bause:2017}, we consider four
test cases based on an L-shaped domain $\Omega=(-0.5,0.5)^2\setminus[0,0.5]^2\subset\mathbb{R}^2$
in the time interval $(0,0.5)$ under two sets of boundary conditions, 
identified by test cases~1a/b/c, and test case~2, cf.~Fig.~\ref{both_mini9:figure:geometry}.
For all cases, vanishing initial conditions $p_0=0$ and $\bm{u}_0=\bm{0}$ are prescribed.
A traction $\bm{\sigma}_\mathrm{N}(t)=(0,-h_\mathrm{max}\cdot256\cdot t^2\cdot(t-0.5)^2)$
is applied on the top with $h_\mathrm{max}$ suitably chosen. Additionally, we prescribe
$p_\mathrm{D}=0$ on the top and $q_\mathrm{N}=0$ on the remaining boundary, zero normal displacement
and homogeneous tangential traction on the left and bottom side, and a homogeneous
traction on the lower right side.
In the test cases~1a/b/c, we also prescribe zero normal displacement on the cut, 
whereas in test case~2, a homogeneous traction is applied on the cut.
The different sets of boundary conditions result in two different physical scenarios.
Despite the two-dimensional, non-symmetric geometry, the test cases~1a/b/c are closely related
to a classical uniaxial compression, whereas the second test case describes a true 
two-dimensional deformation.

For the numerical discretization of Eq.~\eqref{both_mini9:eq:biot:u}--\eqref{both_mini9:eq:biot:q}
in space we use piecewise linear, piecewise constant and lowest order 
Raviart-Thomas finite elements for $\bm{u}$, $p$ and $\bm{q}$, respectively,
defined on a structured, quadrilateral mesh with $16$ elements per $x$- and $y$-direction.
Additionally, for the time discretization we use the backward Euler method
with a fixed time step size $\Delta t=0.01$. 
We solve the discretized Eq.~\eqref{both_mini9:eq:biot:u}--\eqref{both_mini9:eq:biot:q}
using the fixed-stress splitting scheme~\eqref{both_mini9:eq:fs:p}--\eqref{both_mini9:eq:fs:u} for a range of
tuning parameters $K_\mathrm{dr} = \omega K_\mathrm{dr}^{1\text{D}}$ for test cases~1a/b/c
and $K_\mathrm{dr} = \omega K_\mathrm{dr}^{2\text{D}}$ for test case~2,
$\omega\in\{0.5,0.51,...,1.3\}$, and present the accumulated number of iterations
required for convergence of the fixed-stress splitting scheme.
In the following, we denote $K_\mathrm{dr}^\star$ to be the $K_\mathrm{dr}$ yielding minimal number 
of iterations for a single test case.
As stopping criterion, we employ the discrete Euclidean norm $\|\cdot\|_{l^2}$ 
for the algebraic increments between two successive solution vectors for each
of the unknown variables $\bm{u}$, $p$ and $\bm{q}$.
The numerical examples are implemented using the deal.II library 
(and are verified by an implementation using the DUNE library).

\subsection*{Test case 1a -- Effective 1d deformation}
We consider a fixed Young's modulus $E=100$ [GPa] and a varying Poisson's ratio
$\nu\in\{0.01,0.1,0.2,0.3,0.4,0.49\}$.
Moreover, we fix $k = 100$ [mD], $\eta=1$ [cP], $\alpha=0.9$, $M=100$ [GPa], $\bm{g}=\bm{0}$ [m/s$^2$].
On top, we apply the normal force $\bm{\sigma}_\mathrm{N}$ with $h_\mathrm{max}=10$ [GPa].
The number of required fixed stress iterations in relation to the tuning parameter
is displayed in Fig.~\ref{figure:testcase1a}. 
Here, the choice $K_\mathrm{dr}=K_\mathrm{dr}^{1\text{D}}$ is suitable independent of the Poisson's ratio, 
confirming that the problem is essentially driven by uniaxial compression.
\input{ex1_Lshape2d_E001_E100.tex}

\subsection*{Test case 1b -- Soft material}
We modify test case~1a and consider now a softer material with Young's modulus $E=1$ [GPa],
yielding a stronger coupling of the mechanics and flow problem.
On top, we apply the normal force $h_\mathrm{max}=0.1$ [GPa], resulting in a comparable maximal displacement 
of the top boundary. Apart from that, we use the same parameters as in test case~1a.
The number of required fixed stress iterations in relation to the tuning parameter
is displayed in Fig.~\ref{figure:testcase1b}.
Compared to the previous test case, the nature of the mechanical problem becomes more two-dimensional,
indicated by $K_\mathrm{dr}^\star$ lying between $K_\mathrm{dr}^{1\text{D}}$ and $K_\mathrm{dr}^{2\text{D}}$.
More precisely, $K_\mathrm{dr}^\star/K_\mathrm{dr}^{1\text{D}}$ depends on $\nu$.
Only for nearly incompressible materials $K_\mathrm{dr}=K_\mathrm{dr}^{1\text{D}}\approx K_\mathrm{dr}^{2\text{D}}$ is a suitable choice.
Hence, we see that the problem's mechanical character can vary with changing material parameters
but fixed boundary conditions.

\subsection*{Test case 1c -- Influence of flow parameter}
We consider a particular example of test case~1b ($E=1$ [GPa], $\nu=0.01$) for varying permeability 
$k\in\{1e$-$1,1e0,...,1e3\}$ [mD]. Apart from that, we use the same parameters as 
in test case~1b.
The number of required fixed stress iterations in relation to the tuning parameter
is displayed in Fig.~\ref{figure:testcase1c}.
We observe that although all possible mechanical input data is fixed (mechanical boundary conditions 
and material parameters), the optimal tuning parameter $K_\mathrm{dr}^\star$ is in general also dependent on 
material parameters associated with the flow problem. In particular, for decreasing permeability, 
the optimal $K_\mathrm{dr}^\star$ increases towards $K_\mathrm{dr}^{1,\text{D}}$.
\input{Lshape2d_ex2_E001_nu001_ex3_E100b.tex}

\subsection*{Test case 2 -- True 2d deformation}
We consider test case~2 with same parameters as in test case~1a, changing only
from Dirichlet to Neumann boundary conditions on single parts of the boundary.
The number of required fixed stress iterations in relation to the tuning parameter
is displayed in Fig.~\ref{figure:testcase2}.
We observe that the optimal tuning parameter $K_\mathrm{dr}^\star$ is essentially equal
to $K_\mathrm{dr}^{2\text{D}}$, indicating that, the boundary conditions generate a 
true two-dimensional deformation and determine fully the optimal tuning parameter
$K_\mathrm{dr}^\star$.

\section{Conclusion}
In general, the a priori choice of an optimal tuning parameter $K_\mathrm{dr}^\star$ is not trivial
as in the above test cases we have observed:
\begin{itemize}
 \item[$\bullet$] The mathematically motivated tuning parameters $K_\mathrm{dr}^{2\times\lambda}$
 and $K_\mathrm{dr}^{2\times 2\text{D}}$ have in general not in the slightest shown to be optimal,
 which can be confirmed without difficulty for e.g. nearly incompressible materials.
 Instead, the optimal $K_\mathrm{dr}^\star$ has been closer related to the physically 
 motivated parameters $K_\mathrm{dr}^{d^\star\text{D}}$.
 
 \item[$\bullet$] Mechanical boundary conditions are able to determine essentially the physics and define the optimal 
 $K_\mathrm{dr}^\star$, cf.~test case~1a/2, which is consistent with~\cite{both_mini9:kim:2011b}. 
 However, they do not necessarily solely determine the optimal $K_\mathrm{dr}^\star$, cf.~test case~1b. 
 Furthermore, although the fixed-stress approach is based on the inexact inversion of
 the mechanics equation~\eqref{both_mini9:eq:biot:u}, also fluid flow properties can influence $K_\mathrm{dr}^\star$, cf.~test case~1c.
 
 \item[$\bullet$] As expected, for nearly incompressible materials, 
 a suitable tuning parameter is given by $K_\mathrm{dr}=\lambda\approx K_\mathrm{dr}^{d^\star\text{D}}
 \approx \frac{1}{2} K_\mathrm{dr}^{2\times\lambda}\approx \frac{1}{2} K_\mathrm{dr}^{2\times2\text{D}}$, cf.~Section~\ref{both_mini9:Sec:FS}.
 
\end{itemize}

All in all, the optimal tuning parameter $K_\mathrm{dr}^\star$ does not solely depend on the Lam\'e parameters,
but also other physical material parameters, the physical character of the problem and numerical discretization parameters.
Thus, future theoretical analysis of $K_\mathrm{dr}^\star$ should also include the effect of those. 
However, in practice, we expect the dependence of $K_\mathrm{dr}^\star$ on the problem's input data to be complex, and therefore
plan to investigate adaptive techniques for determining a locally defined approximation of $K_\mathrm{dr}^\star$.

\section*{Acknowledgments}  
The research contribution of the second author was partially supported by
E.ON Stipendienfonds (Germany) under the grant T0087 29890 17 while visiting
University of Bergen.

\ifx\undefined\bysame
\newcommand{\bysame}{\leavevmode\hbox to3em{\hrulefill}\,}
\fi

\end{document}